\documentclass[12pt]{article}
\usepackage{latexsym} 
\usepackage{amsmath}
\usepackage{amssymb}
\newtheorem{theorem}{Theorem}[section]
\newtheorem{lemma}[theorem]{Lemma}
\newtheorem{corollary}[theorem]{Corollary}
\newtheorem{proposition}[theorem]{Proposition}
\newtheorem{example}[theorem]{Example}
\newtheorem{remark}[theorem]{Remark}

\def\bit{\begin{itemize}}
\def\eit{\end{itemize}}
\reversemarginpar   
\def\bc{\begin{center}}
\def\ec{\end{center}}
\def\bthm{\begin{theorem}}
\def\ethm{\end{theorem}}
\def\bcor{\begin{corollary}}
\def\ecor{\end{corollary}}
\def\bprop{\begin{proposition}}
\def\eprop{\end{proposition}}
\def\blem{\begin{lemma}}
\def\elem{\end{lemma}}
\def\bex{\begin{example} {\rm }
\def\eex{\end{example} }}
\def\brem{\begin{remark}}
\def\erem{\end{remark}}
\def\prf{\noindent{\bf Proof. }}
\def\bdes{\begin{description}}
\def\edes{\end{description}}
\def\ita{\item[(a)]}
\def\itb{\item[(b)]}
\def\itc{\item[(c)]}
\def\itd{\item[(d)]}

\def\iti{\item[(i)]}
\def\itii{\item[(ii)]}
\def\itiii{\item[(iii)]}
\def\itiv{\item[(iv)]}

\def\beq{\begin{equation}}
\def\eeq{\end{equation}}
\def\ben{\begin{enumerate}}
\def\een{\end{enumerate}}
\def\beqar{\begin{eqnarray}}
\def\eeqar{\end{eqnarray}}
\def\beqarr{\begin{eqnarray*}}
\def\eeqarr{\end{eqnarray*}}

\def\RR{{\mathbb R}}  

\def\qed{\hspace{.1in}{$\blacksquare$} \\}

\def\Pr{{\mathbf{P}}}
\def\E{{\mathbf{E}}}

\title{On strict convergence of stochastic gradients}
 \author{Michel Benaim\\
 Institut de Math\'ematiques\\Universit\'e de Neuch\^atel, Switzerland}

\begin{document}
\maketitle
\begin{abstract}
We discuss conditions ensuring the (strict) convergence of stochastic gradient algorithms.
\end{abstract}


\section{Introduction}
A stochastic gradient algorithm is a process $(x_n), x_n \in \RR^m,$ verifying a recursion of the form
\beq
\label{eq:stochgrad} x_{n+1} - x_n = - \gamma_{n+1} (\nabla V(x_n) + U_{n+1}) \eeq
where $V : \RR^m \mapsto \RR$ is a smooth potential, $(\gamma_n)$ a decreasing sequence of nonnegative weights (e.g $\gamma_n = \frac{A}{n}$) and $(U_n)$ a $\RR^m$ valued sequence of  perturbations (e.g~ a martingale difference sequence).

Under classical assumptions on $(\gamma_n)$ and $(U_n),$ the limit set of $(x_n)$ is a compact connected subset of the {\em critical set} of $V:$
$$\mathsf{crit}(V) = \{x \in \RR^m: \: \nabla V (x) = 0 \}.$$ For generic $V,$ critical points are isolated, so that  the process converges to one of them.
The purpose of this note is to discuss conditions ensuring convergence for general (possibly degenerate) potentials. The main result follows from Lojasiewick type convergence results for gradient like systems established in \cite{Chill09}, combined with a shadowing theorem proved in \cite{Pilyugin99} and classical estimates for stochastic approximation processes.

As an illustration we get assertion $(ii)$ of Theorem \ref{th:illustr} below. Assertion $(i)$ is classical (see e.g~Proposition \ref{th:sard} and Remark \ref{rem:sard} below) and recalled here for completeness.

For $p \in \mathsf{crit}(V)$  let  $D^2 V(p)$ be the Hessian of $V$ at $p$ and let $\Lambda(p)$ denote the spectrum (i.e the collection of eigenvalues) of $- D^2V(p).$ For any set $C \subset \mathsf{crit}(V)$ let $\Lambda(C) = \bigcup_{p \in C} \Lambda(p)$ and ${\cal R}(C) = \RR \setminus \Lambda(C).$
\bthm
\label{th:illustr}
Let $(\Omega, {\cal F}, \Pr)$ be a probability space equipped with a filtration $\{{\cal F}_n\},$  $(U_n)$ a sequence of $\RR^m$ valued random variables adapted to $\{{\cal F}_n\}$ and $x_0 \in \RR^m$ an ${\cal F}_0$ measurable random variable.
Assume
\bdes
\ita  $\gamma_n = \frac{A}{n}$ for some $A > 0;$
 \itb  $\E(U_{n+1} | {\cal F}_n) = 0;$
   \itc $(U_n)$ is bounded in $L^2$ (i.e $sup_n \|U_n\|^2 < \infty$);
   \itd The process $(x_n)$ solution to  (\ref{eq:stochgrad}) is almost surely bounded.
    \edes

Let ${\cal L}(x_n)$ denote the  limit set of $(x_n).$
\bdes
\iti Suppose $V$ is $C^r$ for some $r \geq m.$ Then ${\cal L}(x_n)$ is a compact connected subset of  $\mathsf{crit}(V).$
 \itii Suppose $V$ is real analytic. Let $C$ be a compact connected subset of  $\mathsf{crit}(V).$ Assume that
 \beq
 \label{eq:spec-cond}
 ]-\frac{1}{2A}, 0[ \bigcap {\cal R}(C) \neq \emptyset.
 \eeq Then there exists a random variable $x_{\infty} \in C$ and $c > 0$ such that
 $$\|x_n - x_{\infty}\| = O(\frac{1}{\log(n)^c})$$
 almost surely on the even $$\{{\cal L}(x_n) \subset C\}.$$
 \edes
 \ethm
 \brem
 {\rm If $V$ is proper, there exists $A^* > 0$ such that for all $A \leq A^*$ and every component $C \subset \mathsf{crit}(V)$ condition (\ref{eq:spec-cond}) is satisfied.}
 \erem
 A point $p \in \mathsf{crit}(V)$ is called {\em linearly unstable} if $\Lambda(p) \cap \RR^+_* \neq \emptyset.$ By analogy we say that $C \subset \mathsf{crit}(V)$ is {\em linearly unstable} if every $p \in C$ is linearly unstable. Observe that whenever $C \subset \mathsf{crit}(V)$ is a connected linearly unstable set, then either $C$ reduces to a singleton, or  $0 \in \Lambda(p)$ for all $p \in C.$

 A byproduct of the previous theorem is the following non convergence result.
 \bthm
 \label{th:nonconvgrad}
 With the notation of Theorem \ref{th:illustr}, let $V$ and $C$ be like in Theorem \ref{th:illustr} $(ii).$ Assume furthermore that $C$ is linearly unstable and that
 \bdes \iti $\sup_n E(\|U_n\|^q | {\cal F}_n) < \infty$ for some $q > 2;$
  \itii  $$\liminf_{n \rightarrow \infty} \lambda_{min} [\E(U_{n+1}  U_{n+1}^T | {\cal F}_n)] > 0$$
  where $\lambda_{min}[.]$  denotes  the smallest eigenvalue.
  \edes Then, the event $\{{\cal L}(x_n) \subset C \}$ has zero probability.
 \ethm
\paragraph{Outline} Section \ref{sec:background} introduces some notation and  background. The main result (Theorem \ref{th:main}) is given in section \ref{sec:main} and is applied in  section \ref{sec:nonconv} to provide condition ensuring that the process cannot converge toward an degenerate set of unstable equilibria (Theorem \ref{th:nonconv}).
\section{Notation and background}
\label{sec:background}
Let $F : \RR^m \mapsto \RR^m$ be a $C^1$ globally integrable vector field, $\Phi = \{\Phi_t\}$ the induced flow and $\mathsf{Eq}(F) = \{x \in \RR^m \: : F(x) = 0\}$ the {\em equilibrium set}.

A {\em strict Lyapounov} function for  $\Phi$ (or $F$)  is a continuous map $V : \RR^m \mapsto \RR$ such that
$$V(\Phi_t(x) ) < V(x)$$ for all $x \in \RR^m \setminus \mathsf{Eq}(F)$ and $t > 0.$
When such a $V$ exists, $F$ is called a {\em gradient-like} vector field.

Using the terminology coined in  \cite{BH96},  a continuous function $X : \RR^+ \mapsto \RR^m$ is called a {\em asymptotic pseudo trajectory} (in short APT) of $\Phi,$   if
for all $T > 0$
$$\lim_{t \rightarrow \infty} \sup_{0 \leq h \leq T} \|X(t+h) - \Phi_h(X(t))\| = 0.$$
Its {\em limit set} is the (possibly empty) set defined as
$${\cal L}(X) = \bigcap_{t \geq 0} \overline{ \{X(t+s), \, s \geq 0\}}.$$
\bex
\label{discreteprocess} {\rm
Let $(x_n)$ be solution to the recursion
$$x_{n+1} - x_n = \gamma_{n+1} (F(x_n) + U_{n+1})$$
where $\gamma_n \geq 0, \gamma_n  \rightarrow 0$ and  $\sum_n \gamma_n = \infty.$

Let
$\tau_0 = 0, \tau_n = \sum_{i = 1}^n \gamma_i$ and let $X$ be the continuous interpolated process defined by
\bdes
\ita $X(\tau_n) = x_n,$
\itb $X$ is affine on each interval $[\tau_n, \tau_{n+1}].$
\edes
Assume that
\bdes
\iti $F$ is Lipschitz and bounded on a neighborhood of $(x_n)$  (this holds for instance whenever $(x_n)$ is bounded),
\itii $$\lim_{n \rightarrow \infty} \sup\{ \|\sum_{i = n}^{k-1} \gamma_{i+1} U_{i+1}\|: \: k \geq 1,   \tau_k \leq  \tau_n + T \} = 0.$$
\edes
Then $X$ is an asymptotic trajectory of the flow induced by $F$ (see Proposition 4.1 in \cite{B99}).

}
\eex
\bex[Robbins Monro Algorithm]
\label{ex:RM}
{\rm With the notation of example \ref{discreteprocess}, assume that
\bdes
\iti The sequence $(U_n)$ is a stochastic process adapted to  some filtration $\{{\cal F}_n\},$ $x_0$ is ${\cal F}_0$ measurable and $$\E(U_{n+1} | {\cal F}_n) = 0;$$
\itii $\gamma_n \geq 0,   \sum_n \gamma_n = \infty;$
\itiii There exists $q \geq 2$ such that  $\sum_n \gamma_n^{1 + q/2} < \infty$ and $$\sup_n \E(\|U_n\|^q) < \infty;$$
\itiv $F$ is Lipschitz and bounded on a neighborhood of $(x_n).$
\edes
Then $X$ is an asymptotic trajectory of the flow induced by $F$ (see Proposition 4.3 in \cite{B99}).
}
\eex
By Theorem 0.1 in \cite{BH96} (see also  Theorem 5.7 in \cite{B99}), limit sets of bounded APTs are are {\em internally chain transitive} for $\Phi.$ By this we mean that, if  ${\cal L}$ is  such a set, then  ${\cal L}$ is compact,  connected $\Phi$-invariant  and the restricted flow $\Phi|{\cal L}$ is chain recurrent in the sense of Conley \cite{conley78}.

When $F$ is gradient like this implies  the following result (see Proposition 6.4 in \cite{B99}).

\begin{proposition}
\label{th:sard} Suppose $F$ is gradient like with strict Lyapounov function $V$ and induced flow $\Phi.$ Let $X$ be a bounded APT for $\Phi$ and let ${\cal L} = {\cal L}(X)$ be its
 limit set.
If $V({\cal L} \cap \mathsf{Eq}(F))$ has empty interior, then ${\cal L}$ is a (compact connected) subset of $\mathsf{Eq}(F)$ and $V|{\cal L}$ is constant.

\eprop
\brem
\label{rem:sard}
 {\rm If $\mathsf{Eq}(F) \cap {\cal L}$ is countable (for instance if equilibria of $F$ are isolated) Proposition \ref{th:sard} implies convergence of $(X(t))$ to an equilibrium point.

If $\mathsf{Eq}(F) \subset \mathsf{Crit}(V)$ and $V$ is $C^r $ with $r \geq m.$ Then Sard's theorem implies that $V(\mathsf{Eq}(F))$ has empty interior and the conclusion of Proposition \ref{th:sard} holds.
}\erem

\section{A strict convergence result}
\label{sec:main}
We assume here that $F$ is a $C^1$ vector field which is  gradient like with a $C^1$ strict Lyapounov function $V.$

We let  $X$ denote a
bounded APT.
The {\em error rate} of $X$ is the number $e(X) \in [-\infty, 0]$ defined as
$$e(X) = \limsup_{t \rightarrow \infty} \frac{1}{t}\log (\sup_{0 \leq h \leq T} \|X(t+h) - \Phi_h(X(t))\|)) \leq 0$$
 By Lemma 8.2 in \cite{B99},  $e(X)$ doesn't depend of the choice of  $T > 0.$

If $e(X) \leq \lambda < 0.$ Then, following \cite{Hirsch94} and \cite{BH96}, $X$ is called a {\em $\lambda$ pseudotrajectory}.
\bex[Robbins Monro, continued]
\label{RM3}
{\rm With the notation of example \ref{ex:RM} and $\gamma_n = \frac{A}{n \log(n)^{\beta}}, A >0, 0 \leq \beta \leq 1.$ Then
$$e(X) = - \frac{1}{2A}$$ if $\beta = 0;$ and
$$e(X) = - \infty$$ if $\beta > 0.$ This follows from  Proposition 8.3 and Remark 8.4 in  \cite{B99}.
}
\eex

The map $V$ is said to satisfy a {\em Lojasiewick inequality} at point $p \in \RR^m$ if there exists a neighborhood $U$ of $p$ and constants
$0 < \theta \leq 1/2$ and $c_0 \geq 0$ such that
for all $x \in U$
$$\|V(x) - V(p) \|^{1-\theta} \leq c_0 \|\nabla V(x)\|.$$
It was proved by Lojasiewicz (\cite{Loja}) that $V$ satisfies such an inequality whenever
it is real analytic in a neighborhood of $p$ and this inequality was used to prove that bounded trajectories of real analytic gradient vector fields have finite length, hence converge.
This result was later extended in \cite{Chill09} to gradient-like systems verifying a certain {\em angle condition}.
We say that $(F,V)$ satisfies an angle condition at $p$ if there exists a neighborhood $U$ of $p$ and a constant $c_1 > 0$ such that
$$|\langle \nabla V(x), F(x) \rangle| \geq c_1 \|\nabla V(x)\|\|F(x)\|$$ for all $x \in U.$


Let ${\cal L} = {\cal L}(X)$ denote the limit set of $X$ and let ${\cal R}({\cal L}) \subset \RR$ denote its {\em Sacker-Sell resolvent } (also called in the literature the {\em dynamical resolvent}). Because ${\cal L}$ is internally chain transitive, it follows from Lemma 3 in \cite{Sack78}, that  ${\cal R}({\cal L})$ can be defined as the set of $\lambda \in \RR$
such that for all $x \in {\cal L}$ and $v \in \RR^m \setminus \{0\}$
$$\sup_t \|e^{-\lambda t} D\Phi_t(x) v\| = \infty.$$
In other words, $\lambda \in {\cal R}({\cal L})$ means that there is no bounded solution to the differential system
$$\left \{
  \begin{array}{l}
  \frac{dx}{dt} =  F(x) \\
  \frac{dv}{dt}  =  DF(x) v - \lambda v
  \end{array}
\right. $$
with initial condition $(x,v) \in {\cal L} \times \RR^m \setminus\{0\}.$

For each $p \in \mathsf{Eq}(F)$ let $\Lambda(p)$ denote the set of real parts of eigenvalues of $DF(p)$ (the jacobian matrix of $F$ at $p$). The next proposition can be used to compute or estimate ${\cal R}({\cal L}).$ Its proof follows directly from Proposition \ref{th:sard}, Remark \ref{rem:sard} and the above definition of ${\cal R}({\cal L}).$
\bprop Consider the following assertions:
\bdes \iti $V$ is $C^r$ for $r \geq m$ and $\mathsf{Eq}(F) \subset \mathsf{crit}(V),$
\itii $V({\cal L} \cap \mathsf{Eq}(F))$ has empty interior,
\itiii ${\cal L} \subset \mathsf{Eq}(F),$
\itiv ${\cal R}({\cal L}) = \RR \setminus \bigcup_{p \in {\cal L} \cap \mathsf{Eq}(F)} \Lambda(p).$
\edes
Then $(i) \Rightarrow (ii) \Rightarrow (iii) \Rightarrow (iv).$
\eprop
\bthm
\label{th:main}
Assume that :
\bdes
\iti $V$ satisfies a Lojasiewick inequality and an angle condition at some  point $p \in {\cal L} \cap \mathsf{Eq}(F);$
\itii $]e(X), 0[ \cap {\cal R}({\cal L}) \neq \emptyset.$
\edes
Then ${\cal L}$ reduces to $\{p\}.$
If furthermore $|\langle \nabla V(x), F(x) \rangle| \geq \beta \|\nabla V(x)\|^2$ in a neighborhood of $p$, then
$$\|X(t) - p\| = \left \{ \begin{array}{c}
                            O(e^{ct}) \mbox { if } \theta = 1/2 \\
                            O(t^{-\theta/(1-2 \theta)}) \mbox { if } 0 < \theta < 1/2
                          \end{array} \right.$$
                          where $c = \max(-\frac{\beta}{2c_0^2}, \inf ]e(X),0[ \cap {\cal R}({\cal L})).$
\ethm

\prf
For every $\RR^m$-valued sequence $\xi = (\xi_k)_{k \geq 0}$ and every $r \geq 1$ set
$$\|\xi\|_r = \sum_{k \geq 0} r^k \|y_k\|.$$
Given such a sequence and $x \in \RR^m,$ define the sequences $$g(\xi)= (g_k(\xi))$$  and $$h(x,\xi) = (h_k(x,\xi))$$ by
$$g_k(\xi) = \xi_{k+1} - \Phi_1(\xi_k)$$ and
$$h_k(x,\xi) = \Phi_k(x) - \xi_k$$ for all $k \geq 0.$

Let $C$ be a compact connected invariant set (say internally chain transitive),  $\mu  \in {\cal R}(C) \cap \RR^-$ and $r = e^{-\mu}.$ The following result follows from the more general
 Theorem 1.4.5 in \cite{Pilyugin99}
 \bprop[Theorem 1.4.5 in \cite{Pilyugin99}]  There exist positive numbers $d_0, L$ and a neighborhood $W$ of $C$ such that if the sequence $\xi$ is contained in $W$ and verifies
$\|g(\xi)\|_r \leq d \leq d_0$ then there exists $x \in \RR^m$ such that  $\|h(x,\xi)\|_r \leq L d.$
\eprop
In order to use this proposition, let $$\xi_k = X(T + k)$$ where $T > 0$ will be chosen later and let $C = {\cal L}$ be the limit set of $X.$

Choose (thank to condition $(ii)$)  $\mu \in ]e(X), 0[ \cap {\cal R}(C)$ and $e(X) < \alpha < \mu.$ Then, for $t$ large enough
$$\sup_{0 \leq h \leq 1} \|X(t+h) - \Phi_h(X(t))\| \leq e^{\alpha t}.$$ Therefore, with $r = e^{-\mu}$ and $T$ sufficiently large,
$\|g(\xi)\|_r \leq \frac{e^{\alpha T}}{1 - e^{\alpha - \mu}} \leq d.$  Thus, by the latter proposition,  $\|h(x,\xi)\|_r \leq L d$ for some $x.$ It then follows that for all $0 \leq h \leq 1$ and $k \geq 0$
$$\|X(T+ k + h) - \Phi_{k+h}(x)\| \leq \|X(T+ k + h) - \Phi_h(X(T+k))\| + \|\Phi_h(X(T+k)) - \Phi_h (\Phi_k(x))\|$$
$$\leq e^{\alpha (T+k)} + e^L \|g_k(\xi)\| = o(e^{\mu k})$$ where $L$ is a Lipschitz constant of $F$ on a neighborhood of $X(\RR^+).$ Then
$$\|X(T+t) - \Phi_t(x)\| = o(e^{\mu t}).$$
In particular $X(.)$ and $\{\Phi_t(x), t \geq 0 \}$ have the same limit set.
Now, by assumption $(i)$ and Theorem 1 in \cite{Chill09}\footnote{This theorem is stated under the assumption that $V$ is $C^2$ and satisfies a global angle condition, but the proof works for $C^1$ with an angle condition at point $p$}, $\Phi_t(x) \rightarrow p$ as $t \rightarrow \infty.$ Under the supplementary assumption that $|\langle \nabla V(x), F(x) \rangle| \geq \beta \|\nabla V(x)\|^2$ then, by Theorem 2 in \cite{Chill09}
$$\|\Phi_t(x) - p\| = \left \{ \begin{array}{c}
                            O(e^{-ct}) \mbox { if } \theta = 1/2 \\
                            O(t^{-\theta/(1-2 \theta)} \mbox { if } 0 < \theta < 1/2
                          \end{array} \right.$$

\qed
\paragraph{Proof of Theorem  \ref{th:illustr}} Assertion $(i)$ follows from Remark \ref{rem:sard} and Proposition \ref{th:sard}. Assertion $(ii)$
follows from Theorem \ref{th:main}, Example \ref{RM3} and the fact that real analytic maps satisfies Lojasiewick inequality. \qed
\paragraph{Discussion of condition (ii) in Theorem \ref{th:main}}
If we no longer assume that $F$ is gradient like (but continue to assume that $X$ is a bounded APT with limit set ${\cal L} = {\cal L}(X)$)
  the condition
 \beq
 \label{eq:gapSack}
 ]e(X), 0[ \cap {\cal R}({\cal L}) \neq \emptyset
 \eeq
 implies (by Theorem 1.4.5 in \cite{Pilyugin99} as in the proof of Theorem \ref{th:main}) that there exists $x \in \RR^m$ such that
 $$\limsup_{t \rightarrow \infty} \frac{\log( \|X(t) - \Phi_t(x)\|}{t} \leq \lambda$$
 with $\lambda = \inf  ]e(X), 0[ \cap {\cal R}({\cal L}) < 0.$

 This shadowing property is reminiscent of  the shadowing Theorem 9.3 in \cite{BH96} (see also Theorem 8.9 in \cite{B99}). This latter result which easily follows from   Hirsch's shadowing theorem (\cite{Hirsch94}, Theorem 3.2) asserts the following.
Let $${\cal E}(\Phi, {\cal L}) = \lim_{t \rightarrow \infty} \frac{1}{t} \log (\inf_{x \in {\cal L}} \| D\Phi_{-t} (\Phi_t(x)) \|^{-1})$$ be the {\em expansion rate} of $\Phi$ at ${\cal L}.$
If
 \beq
 \label{eq:gapHirsch}
 e(X)  < \min(0, {\cal E}(\Phi, {\cal L}))
 \eeq then there exists $x \in \RR^m$ such that
 $$\limsup_{t \rightarrow \infty} \frac{1}{t} \log( \|X(t) - \Phi_t(x) \|) \leq e(X).$$
By a theorem of Schreiber (\cite{Schreiber97})
$${\cal E}(\Phi, {\cal L}) = \inf_{\mu \in {\cal P}_{erg}({\cal L})} \lambda_1(\mu) $$
where ${\cal P}_{erg}({\cal L})$ is the set of $\Phi$-invariant ergodic measures supported by ${\cal L}$ and for each
 $\mu \in {\cal P}_{erg}({\cal L})$ $\lambda_1(\mu)$ is the smallest Lyapounov exponent
\footnote{By a lyapounov exponent of $\mu$ we mean  a lyapounov exponent of the skew product
flow $(x,v) \in {\cal L} \times \RR^m \mapsto (\Phi_t(x), D\Phi_t(x) v)$}
of $\mu.$

Now, by Theorem 2 in \cite{Sack78}
the {\em dynamical spectrum}  $\Sigma({\cal L}) = \RR \setminus {\cal R}({\cal L})$ is the union of $k \leq m$ compact intervals
$$\Sigma({\cal L}) = [a_1,b_1] \cup \ldots \cup [a_k,b_k]$$
with $a_1 \leq b_1 < a_2 \leq b_2 < \ldots < a_k \leq b_k;$ and
by Theorem 2.3 in \cite{Johnson87} for every $\mu \in {\cal P}_{erg}({\cal L})$ all the lyapounov exponents of  $\mu$ are contained in $\Sigma({\cal L})$ and  every
point in $\partial \Sigma({\cal L}) = \{a_1,b_1, \ldots, a_k,b_k\}$ is a  lyapounov exponent for some $\mu \in {\cal P}_{erg}({\cal L}).$
It then follows that
$$a_1 = {\cal E}(\Phi, {\cal L})$$ so that condition (\ref{eq:gapHirsch}) implies condition (\ref{eq:gapSack}).
\section{Non convergence toward (degenerate) unstable set of equilibria}
\label{sec:nonconv}
We discuss here some implication of Theorem \ref{th:main} to the problem of non convergence toward degenerate set of equilibria, for stochastic approximation processes.

Let $X = (X(t))_{t \geq 0}$ be a $\RR^m$ valued continuous  stochastic process. For every open set $U \subset \RR^m$ and $s \geq  0,$ let
 $T_s^U$ be the stopping time (with respect to the canonical filtration of $X$) defined by
 $$T_s^U = \inf \{ t \geq s \: : X(t) \in \RR^m \setminus U \}.$$ Let $p \in \RR^m.$
 We say that $p$ is {\em repulsive} for $X$ provided there exists a  neighborhood $U$ of $p$ such that for all $n \in \mathbb{N}$
  $$\Pr( T_n^U < \infty) = 1.$$
 This clearly implies that
$$\Pr(\lim_{t \rightarrow \infty} X(t) = p) = 0$$
but the converse is  false as shown by the next example.
\bex {\rm
Consider a two colors Polya urn process. Initially, there are two balls, say one black and one white, in the urn. At each time a ball is randomly chosen, and replaced in the urn with a new ball of the same color. Let $W_n$ denote the number of white balls at time $n$ and $x_n = \frac{W_n}{n+2} \in [0,1]$ its proportion. Then $(x_n)$ satisfies (\ref{discreteprocess})
with $\gamma_n = \frac{1}{n+2}$ and $F(x) = 0.$ It is well known that $(x_n)$ (hence the interpolated process $X$) converges almost surely toward a random variable  $X_{\infty}$ having a uniform distribution on $[0,1]$. Thus, for all $p \in [0,1]$ $\Pr(\lim_{t \rightarrow \infty} X(t) = p) = 0$ but $p$ is not repulsive for $X$.}
\eex
  The following result is a straightforward consequence of Theorem \ref{th:main}.
  \bthm
\label{th:nonconv}
Let $F$ and $V$ be like in section \ref{sec:main} and let  $X : \RR^+ \mapsto \RR^m$  be a continuous stochastic process, almost surely bounded and such that $e(X) \leq e < 0.$ Let $C \subset \mathsf{Eq}(F)$ be a compact connected set of equilibria. Assume that
\bdes
\iti $V$ satisfies a Lojasiewick inequality and an angle condition at every point $p \in C;$
\itii $]e, 0[ \cap {\cal R}(C)) \neq \emptyset;$
\itiii Every point $p \in C$ is repulsive for $X.$
\edes
Then $$\Pr({\cal L}(X) \subset C) = 0.$$
\ethm
\prf By assumptions and compactness, there exist points $p_1, \ldots, p_k \in C$ and neighborhoods $U_1, \ldots, U_k$ such that $C \subset \bigcup U_i, p_i \in U_i,$ and $\Pr({\cal L}(X) \subset U_i) = 0.$ Now, on the event ${\cal L}(X) \subset C$ there exists $i \in \{1, \ldots, k\}$ such that ${\cal L}(X) \cap U_i \neq \emptyset,$ hence  by Theorem \ref{th:main}, ${\cal L}(X) \subset U_i.$ Thus $\Pr({\cal L}(X) \subset C) \leq \sum_i \Pr( {\cal L}(X) \subset U_i) = 0.$
\qed

Theorem \ref{th:nonconv} requires the verification that every point in $C$ is repulsive.
Beginning with a seminal paper by Pemantle \cite{Pem90}, the literature on stochastic approximation and urn processes has produced several results showing  that, under reasonable assumptions, a process given by (\ref{discreteprocess}) cannot converge toward an "unstable" equilibrium of the associated vector field. It turns out that, while these results are usually formulated as "non convergence results", a careful reading of the proofs shows that they actually prove repulsiveness of the unstable point.

The following result, due to Brandi\`ere and Duflo \cite{BraDuf96} (see also Chapter 3 in \cite{Duf96})  was first proved by Pemantle \cite{Pem90} when $F$ is $C^2$ and $p$  hyperbolic, for bounded noise and  $\gamma_n = 1/n.$ Pemantle's theorem was later extended to more general invariant sets (including non hyperbolic equilibria) in \cite{B99}, (Section 9, Theorem 9.1). Under stronger assumptions on the noise sequence the smoothness assumption on $F$ can be weakened to $C^1$ (see Theorem 3.12 in \cite{BF12}). Tarr\`es' PhD thesis  (\cite{Tar}) contains  interesting generalizations that allow to deal with  unstable (but not linearly unstable) points, especially in dimension one.

Let $p$ be an equilibrium of  a $C^1$ vector field $F : \RR^m \mapsto \RR^m.$  Then $\RR^m$ can be written as the direct sum of $E^s_p, E^u_p, E^c_p$
the generalized eigenspaces corresponding to the eigenvalues of the jacobian matrix
$DF(p)$ having, respectively, negative real parts, null real parts and positive
real parts. Equilibrium $p$ is said to be {\em linearly unstable} if $E^u_p \neq \{0\}$ and $hyperbolic$ if $E^c_p = \{0\}.$
\bthm[\cite{BraDuf96}, Theorem 1]
\label{th:BD}
Let $F$ be a $C^1$ vector field whose Jacobian is  locally Lipschitz. Let $p$ be a linearly unstable (non necessarily hyperbolic) equilibrium of $F.$
Consider the process given in example \ref{discreteprocess} where $(U_n)$ and $\gamma_n$ are  like in example \ref{ex:RM} with $q = 2.$  Assume furthermore  that
$$\sup_n \E (\|U_{n+1}\|^2 | {\cal F}_n) < \infty$$ and
$$\liminf_{ n \rightarrow \infty} \E (\|U^u_{n+1}\| | {\cal F}_n) > 0,$$ where $U_n^u$ stands for the projection of $U_n$ on $E^u_p$ along $E^s_p \oplus E^c_p.$ Then $p$ is repulsive for $X.$
\ethm
\prf  Theorem 1 in  \cite{BraDuf96} states that $\Pr( \lim_{n \rightarrow } x_n = p) = 0$ so we need to explain how the proof implies that $p$ is repulsive. We can always assume $p = 0.$  By the center-stable manifold theorem (see e.g~\cite{Rob99}, Section 5.10.2), there exists a neighborhood $U = U_1 \times U_2$ of the origin with $U_1 \subset E^s_p \oplus E^c_p, U_2 \subset E^u_p,$ and  a $C^1$ map $G : U_1 \mapsto U_2$ with $G(0) = 0, DG(0) = 0$ whose graph $W^{cs}_{loc}$ is locally invariant under $\{\Phi_t\}.$ Furthermore, if $\Phi_t(x) \in U$ for all $t \in \RR,$ then $x \in W^{cs}_{loc}.$
In particular, if $X$ is an asymptotic pseudo trajectory of $\Phi,$ then $${\cal L}(X) \subset U \Rightarrow {\cal L}(X) \subset W^{cs}_{loc}.$$
Let $X$ be the interpolated process associated to  $(x_n)$ (as defined in example \ref{discreteprocess}). Set $X(t)  = (X^1(t), X^2(t)) \in (E^s_p \oplus E^c_p) \times E^u_p.$ The proof of Theorem 1 in  \cite{BraDuf96} shows that
$$\Pr(\lim_{t \rightarrow \infty}  X^2(t) - G(X^1(t)) = 0) = 0.$$
Thus, $$\Pr({\cal L}(X) \subset U) = 0.$$ \qed

\paragraph{Proof of Theorem \ref{th:nonconvgrad}} follows from Theorems \ref{th:nonconv} and \ref{th:BD} and the discussion in (\cite{BraDuf96}, section I.2, 1) which shows that the assumptions on the noise given in Theorem \ref{th:nonconvgrad} imply the assumptions given in Theorem \ref{th:BD}.

\bibliographystyle{amsplain}
\bibliographystyle{imsart-nameyear}
\bibliographystyle{nonumber}
\bibliography{stochgrad}

\providecommand{\bysame}{\leavevmode\hbox to3em{\hrulefill}\thinspace}
\providecommand{\MR}{\relax\ifhmode\unskip\space\fi MR }
\providecommand{\MRhref}[2]{%
  \href{http://www.ams.org/mathscinet-getitem?mr=#1}{#2}
}
\providecommand{\href}[2]{#2}
\begin{thebibliography}{10}

\bibitem{B99}
M~Bena{\"{\i}}m, \emph{Dynamics of stochastic approximation algorithms},
  S\'eminaire de {P}robabilit\'es, {XXXIII}, Lecture Notes in Math., vol. 1709,
  Springer, Berlin, 1999, pp.~1--68.

\bibitem{BF12}
M.~Bena\"{\i}m and M.~Faure, \emph{Stochastic approximations, cooperative
  dynamics and supermodular games}, The Annals of Applied Probability
  \textbf{22} (2012), no.~5, 2133--2164.

\bibitem{BH96}
M~Bena{\"{\i}}m and M~W. Hirsch, \emph{Asymptotic pseudotrajectories and chain
  recurrent flows, with applications}, J. Dynam. Differential Equations
  \textbf{8} (1996), no.~1, 141--176. \MR{1388167 (97d:58165)}

\bibitem{BraDuf96}
O.~Brandiere and M.~Duflo, \emph{{Les algorithmes stochastiques contournent-ils
  les pi{\`e}ges?}}, Annales de l'I. H. P. Probabilit{\'e}s et statistiques
  \textbf{32} (1996), no.~3, 395--427.

\bibitem{Chill09}
R.~Chill, A.~Haraux, and M.~Ali-Jendoubi, \emph{Applications of the lojasiewicz
  simon, gradient inequality to gradient-like evolution equations}, Analysis
  and Applications \textbf{7} (2009), no.~4, 351--372.

\bibitem{conley78}
C.~Conley, \emph{Isolated invariant sets and the morse index}, AMS, Providence,
  1978, Series in Math.

\bibitem{Duf96}
M.~Duflo, \emph{{Algorithmes stochastiques}}, Springer Paris, 1996.

\bibitem{Hirsch94}
M.~W Hirsch, \emph{Asymptotic phase, shadowing and reaction diffusion systems},
  Differential equations, dynamical systems and control science, Lecture Notes
  in Pure and Applied Math., vol. 152, Marcel Dekker, New York, 1994,
  pp.~87--99.

\bibitem{Loja}
S.~Lojasiewicz, \emph{Une propri\'et\'e topologique des sous-ensembles
  analytiques r\'eels}, Les \'Equations aux D\'eriv\'ees Partielles (1963),
  87–--89, \'Editions du C.N.R.S, Paris.

\bibitem{Pem90}
R.~Pemantle, \emph{{Nonconvergence to unstable points in urn models and
  stochastic approximations.}}, ANN. PROB. \textbf{18} (1990), no.~2, 698--712.

\bibitem{Pilyugin99}
S.~Y. Pilyugin, \emph{Shadowing in dynamical systems}, Springer, 1999, Lecture
  Notes in Mathematics, 1706.

\bibitem{Johnson87}
K.~J.~Palmer R.~A.~Johnson and G.~R. Sell, \emph{Eergodic properties of linear
  dynamical systems}, SIAM J. MATH. ANAL. \textbf{18} (1987), no.~1, 1--33.

\bibitem{Rob99}
Clark. Robinson, \emph{{Dynamical Systems, Stability, Symbolic Dynamics, and
  Chaos}}, 1999, second edition.

\bibitem{Sack78}
R.~J. Sacker and G.~R. Sell, \emph{A spectral theory for linear differential
  systems}, Journal of Differential Equations \textbf{27} (1978), 320--358.

\bibitem{Schreiber97}
S.~Schreiber, \emph{Expansion rates and lyapounov exponents}, Discrete and
  Continuous Dynamical Systems \textbf{3} (1997), no.~4, 433--438.

\bibitem{Tar}
P.~Tarr{\`e}s, \emph{{Pi{\`e}ges des algorithmes stochastiques et marches
  al{\'e}atoires renforc{\'e}es par sommets}}, Ph.D. thesis, ENS Cachan, 2001.

\end{thebibliography}

\section*{Acknowledgments}
This work was supported by the SNF grant $2000020 - 149871/1$

\end{document}